\documentclass[12pt]{article}

\usepackage{fullpage}

\usepackage{amsfonts}
\usepackage{graphicx}
\usepackage{cite}
\usepackage{mathrsfs}
\usepackage{amsmath}
\usepackage{amssymb}
\usepackage{float}
\usepackage{color}
\usepackage{setspace}
\usepackage{graphicx}
\usepackage{subfigure}
\usepackage{amsmath}\numberwithin{equation}{section}
\makeatletter
\renewcommand{\@thesubfigure}{\hskip\subfiglabelskip}
\makeatother

\newtheorem{thm}{Theorem}[section]

\newtheorem{remark}[thm]{Remark}

\definecolor{lw}{RGB}{0,0,255}
\definecolor{zs}{RGB}{255,0,0}

\definecolor{lightblue}{RGB}{0.68,0.85,0.9}
\definecolor{bluegray}{rgb}{0.4, 0.6, 0.8}

\def\RW{{\mathrm{RW}}}

\begin{document}


\begin{center}
\Large
{\bf Spectral Radius  of Biased Random Walks \\ on Regular Trees}

\end{center}

\begin{center}
Song He
\vskip 1mm
\footnotesize{School of Mathematics and Statistics, Huaiyin Normal University}\\
\footnotesize{Huaian City, 223300, P. R. China}\\
\footnotesize{~~~~~~~~~~songhe@hytc.edu.cn}\\
\end{center}

\begin{abstract}
We consider biased random walk on regular tree and we obtain the spectral radius, first return probability and $n$-step transition probability.
\end{abstract}

\vspace{0.4cm}

\section{Introduction}
Let $G=(V(G), \, E(G))$ be a locally finite, connected infinite graph, where $V(G)$ is the set of its vertices and
$E(G)$ is the set of its edges. Fix a vertex $o$ of $G$ as {the} root, we assume that $o$ has at least one neighbor.
For any vertex $x$ of $G$ let $|x|$ denote the graph distance between $x$ and $o$. Let
$\mathbb{N} := \{ 1, \, 2, \ldots\}$ and  $\mathbb{Z}_{+}=\mathbb{N}\cup\{0\}$. For any $n\in\mathbb{Z}_{+}$:
$$
B_G(n)
=
\{x\in V(G):\ |x| \le n\},
\qquad
\partial B_G(n)
=
\{x\in V(G):\ |x| =n\}.
$$
Let $d\in\mathbb{N}$, $d\ge2$, $\mathbb{T}_d$ denotes $d$-regular trees. Fix a vertex $o$ of $\mathbb{T}_d$ as {the} root.
For $\lambda>0$, if an edge $e=\{x,y\}$ is at distance $n$ from $o$, its conductance is defined as $\lambda^{-n}$. Denote by $\RW_\lambda$  the nearest-neighbour random walk $(X_n)_{n=0}^\infty$ among such conductances and call it the $\lambda$-{\it biased random walk}. In other words, $\RW_{\lambda}$ on $\mathbb{T}_d$
has the following transition probabilities: for $v \sim u$ (i.e., if $u$ and $v$ are adjacent on $\mathbb{T}_d$),
\begin{eqnarray}
    \label{(1.1)}
    p(v,u)
    :=
    p_{\lambda}^{\mathbb{T}_d}
  =
  \begin{cases}
    \frac{1}{d}
        &\mathrm{if}\ v=o,
        \\
        \frac{\lambda}{d+\lambda-1}
        &\mathrm{if}\ u\in \partial B_{\mathbb{T}_d}(|v|-1)\ \mathrm{and}\ v\neq o,
        \\
        \frac{1}{d+\lambda-1}
        &\mathrm{otherwise}.
  \end{cases}
\end{eqnarray}
Notice $\RW_1$ is just the simple random walk on $\mathbb{T}_d.$  By Rayleigh's monotonicity principle (see \cite{LR-PY2016}, p. 35), there is a critical parameter $\lambda_c(G)\in (0, \, \infty]$ such that $\RW_{\lambda}$ is transient for $\lambda<\lambda_c(G)$ and is recurrent for $\lambda>\lambda_c(G)$.
\medskip
In the following we will introduce some basic notations. Write
$$
p^{(n)}(x,\, y)
:=
p^{(n)}_\lambda(x,\, y)=\mathbb{P}_x(X_n=y),
$$

\noindent where $\mathbb{P}_x:=\mathbb{P}_x^{G}$ is the law of $\RW_\lambda$ starting at $x$. The Green function is given by
\[
\mathbb{G}(x,\, y\, |\, z)
:=
\mathbb{G}_\lambda(x,\, y\, | \, z)
=
\sum_{n=0}^\infty p^{(n)}(x,\, y)z^n,~ x,~y\in V(G),~z\in \mathbb{C},\ |z| < R_{\mathbb{G}} \, ,
\]

\noindent where $R_\mathbb{G}=R_\mathbb{G}(\lambda)=R_\mathbb{G}(\lambda,x,y)$ is its convergence radius. Recall \cite{LR-PY2016} {\it Exercise 1.2},
\[
R_\mathbb{G}
=
R_\mathbb{G}(\lambda)
=
\frac{1}{\limsup_{n\to \infty}\sqrt[n]{p^{(n)}(x,y)}}
\]

\noindent  is independent of $x$, $y$ when $\RW_\lambda$ is irreducible, i.e., $\lambda>0$. Call
\[
\rho_\lambda
=
\rho(\lambda)
=
\frac{1}{R_\mathbb{G}}
=
\limsup_{n\to \infty} p^{(n)}(x,x)^{1/n}
=
\limsup_{n\to \infty} p^{(n)}(o,o)^{1/n}
\]

\noindent the spectral radius of $\RW_\lambda$.
Let $M_n=|\partial B_G(n)|$ be the cardinality of $\partial B_G(n)$ for any $n\in\mathbb{Z}_{+}.$ Define the growth rate of $G$ as
\[\mathrm{gr}(G)=\liminf_{n\rightarrow \infty}\sqrt[n]{M_n}.\]
Since the sequence $\{M_n\}_{n=0}^\infty$ is submultiplicative, the limit $\mathrm{gr}(G)=\lim\limits_{n\rightarrow\infty}\sqrt[n]{M_n}$ exists indeed.

\medskip
The motivation for introducing $\RW_\lambda$ was to design a Monte-Carlo algorithm for self-avoiding walks by Berretti and Sokal \cite{BA-SA-1985}. See \cite{LG-SA1988,SA-JM1989, RD1994} for refinements of this idea. Due to interesting phenomenology and similarities to concrete physical systems (\cite{BM-DD1983, DD1984, DD-SD1998, HS-BA2002, SSSWX2021}), biased random walks and biased diffusions in disordered media have attracted much attention in mathematical and physics communities since the 1980s.

\medskip
In the 1990s, Lyons (\cite{LR1990, LR1992, LR1995}), and Lyons, Pemantle and Peres (\cite{LR-PR-PY1996a, LR-PR-PY1996b}) made series of achievements in the study of $\RW_\lambda$'s. $\RW_\lambda$ has also received attention recently, see \cite{BG-HY-OS2013, AE2014, BG-FA-SV2014, HY-SZ2015} and the references therein. Ben Arous and Fribergh publish a survey on biased random walks on random graphs see\cite{BG-FA2014}. For spectral radius, R. Lyons \cite{LR1992} showed that the critical parameter for
$RW_{\lambda}$ on a general tree is exactly the exponential of the Hausdorff dimension
of the tree boundary. And R. Lyons \cite{LR1995} proved that for Cayley graphs and degree bounded transitive graphs,
the growth rate is exactly the critical parameter of the $RW_{\lambda}.$
This paper focuses on a specific properties of spectral radius of $\RW_\lambda$'s on non-random infinite graphs.

We are ready to state our main results. The proofs will be presented in Section 2.
\begin{thm}\label{thm1.1}
For the $d$-regular tree $\mathbb{T}_d$, the following holds:
$$
\rho_{\mathbb{T}_d}(\lambda)
=
\frac{2\sqrt{(d-1)\lambda}}{d-1+\lambda},
\qquad
\lambda\in (0, \, \lambda_c(\mathbb{T}_d)]=(0,\, d-1],
$$
and for $\lambda\in (0, \, \infty)$ and $n\to \infty$,
\begin{equation}
    f^{(2n)}_\lambda(o, \, o)
    \sim
    \frac{1}{\sqrt{\pi}} \left( \frac{2\sqrt{(d-1)\lambda}}{d-1+\lambda}\right)^{\!\! 2n} n^{-3/2}.
    \label{f_regular_tree}
\end{equation}
Moreover,
\begin{equation}
    p^{(2n)}_\lambda(o,o)
    \sim
    \begin{cases}
          \frac{(d-1-\lambda)^2}{16(\pi\lambda)^{1/2}(d-1)^{3/2}}\rho_{\mathbb{T}_d}(\lambda)^{2n}n^{-3/2}&\hbox{if}\ \lambda\in (0, \, d-1),
      \\
    \frac{1}{\sqrt{\pi n}}&\hbox{if}\ \lambda=d-1.
    \end{cases}
    \label{p_regular_tree}
\end{equation}
\end{thm}
\begin{remark}
Since for the case $\lambda>\lambda_c(\mathbb{T}_d)=d-1$,  $\RW_\lambda$ is recurrence, it means that $\rho_\lambda=1$. Hence, the spectral radius $\rho_\lambda$ is continuous in $\lambda\in
(0,\infty)$, {and}
$\rho(\lambda_c(\mathbb{T}_d)) = 1$.
\end{remark}

\begin{remark}
For $\lambda\in (0,\, d-1)$, the derivative of $\rho_{\mathbb{T}_d}(\lambda)$:
$$
\rho^{'}_{\mathbb{T}_d}(\lambda)=\frac{\sqrt{(d-1)/\lambda}(d-1+\lambda)}{(d-1-\lambda)^2}>0.
$$
It means that $\rho_{\mathbb{T}_d}(\lambda)$ is strictly increasing for $\lambda\in (0,\, d-1)$.
\end{remark}

\section{Proof of Therom~\ref{thm1.1}}\label{s:pfthm2.2}

\noindent\textbf{Proof of Theorem~\ref{thm1.1}} Assume $\lambda>0.$ Notice that $\RW_\lambda$ $(X_n)_{n=0}^{\infty}$ must return to $o$ in even steps, and that $\{|X_n| \}_{n=0}^{\infty}$ with $|X_0| =0$
is a Markov chain on $\mathbb{Z}_{+}$ with transition probabilities given by
$$
p(x, \, y)
=
\left\{\begin{array}{cl}
1 &{\rm if}\ x=0, \ y=1\\
      \frac{\lambda}{d-1+\lambda} &{\rm if}\ y=x-1\ \mbox{and}\ x\neq 0,\\
     \frac{d-1}{d-1+\lambda} &{\rm otherwise}.
   \end{array}
\right.
$$
\noindent Recall for any $n\in \mathbb{N}$ and $k\in\mathbb{Z}_{+}$,
$$
f^{(2n)}_\lambda (o,\, o)
=
\mathbb{P}_o\left(\tau_o^{+}=2n\right),
\qquad
f^{(2n-1)}_\lambda(o,\, o)=0,\ \lambda\in (0,\, \infty),
$$

\noindent and the $k$th Catalan number given by $c_k = \frac{1}{k+1}{{2k}\choose{k}}$, with the associated related generating function
\begin{equation}
    \mathcal{C}(x)
    :=
    \sum_{k=0}^{\infty}c_kx^k
    =
    \frac{1-\sqrt{1-4x}}{2x},
    \qquad
    x\in \left[-\frac{1}{4}, \, \frac{1}{4}\right].
    \label{Catalan_generating_fct}
\end{equation}

\noindent Note the number of all $2n$-length nearest-neighbour paths $\gamma=w_0w_1\cdots w_{2n}$ on $\mathbb{Z}_{+}$ such that
$$
w_0=w_{2n}=0,
\qquad
w_j\ge 1,\ 1\le j\le 2n-1
$$

\noindent is precisely $c_{n-1}$. Hence for any $\lambda>0,$
$$
f_\lambda^{(2n)}(o, \, o)
=
c_{n-1}\left(\frac{d-1}{d-1+\lambda}\right)^{n-1}\left(\frac{\lambda}{d-1+\lambda}\right)^n,\ n\in\mathbb{N},
$$

\noindent which readily yields \eqref{f_regular_tree} by means of Stirling's formula.

By definition, for $\lambda>0$,
\begin{align*}
    U_\lambda(o,\, o\, |\, z)
 &=\sum_{n=1}^\infty f^{(2n)}_\lambda(o,\, o)z^{2n}
    =
    \sum_{n=1}^\infty c_{n-1}
    \left(\frac{d-1}{d-1+\lambda}\right)^{n-1}
    \left(\frac{\lambda}{d-1+\lambda}\right)^n z^{2n}
    \\
 &=\frac{\lambda}{d-1+\lambda} z^2 \, \mathcal{C}\left(\frac{\lambda(d-1)z^2}{(d-1+\lambda)^2}\right) ,
\end{align*}

\noindent which, in view of \eqref{Catalan_generating_fct}, implies that for $|z| \le \frac{d-1+\lambda}{2\sqrt{\lambda(d-1)}}$,
\begin{equation}
    U_\lambda(o,\, o\, |\, z)
    =
    \frac{(d-1+\lambda)-\sqrt{(d-1+\lambda)^2-4\lambda(d-1)z^2}}{2(d-1)} .
    \label{U_regular_tree}
\end{equation}
\noindent Taking $z=1$ gives that
$$
\mathbb{P}_o\left(\tau_o^{+}<\infty \right)
=
U_\lambda(o, \, o \, | \, 1)
=
\frac{\lambda\wedge (d-1)}{d-1}.
$$
Notice from \eqref{U_regular_tree} that when $0<\lambda\le d-1$,
$$
U_{\lambda}\left(o,\, o \, \Big| \, \frac{d-1+\lambda}{2\sqrt{\lambda(d-1)}} \right)
=
\frac{d-1+\lambda}{2(d-1)}
\le 1.
$$

\noindent Hence, for $|z|<\frac{d-1+\lambda}{2\sqrt{\lambda(d-1)}}$ and $0<\lambda\le d-1$,
\begin{align}
    \mathbb{G}_\lambda(o, \, o\, |\, z)
 &=\frac{1}{1-U_\lambda(o, \, o\, |\, z)}
    \nonumber
    \\
 &=\frac{2(d-1)}{2(d-1)- (d-1+\lambda)+\sqrt{(d-1+\lambda)^2-4\lambda(d-1)z^2}} .
    \label{G_regular_tree}
\end{align}
\noindent This implies that the convergence radius for $\mathbb{G}_\lambda(o, \, o\, |\, z)$ is $\frac{d-1+\lambda}{2\sqrt{\lambda(d-1)}}$. In other words,
$$
\rho(\lambda)
:=
\rho_{\mathbb{T}_d}(\lambda)
=
\frac{2\sqrt{\lambda(d-1)}}{d-1+\lambda},
\qquad
0<\lambda\le d-1.
$$

It remains to show \eqref{p_regular_tree} for $\lambda\in (0, \, d-1)$. Write $a(\lambda)=\frac{2(d-1)}{d-1+\lambda}$ and $b(\lambda)=\frac{d-1-\lambda}{d-1+\lambda}.$ Then for any $|z|\le R_{\mathbb{G}}(\lambda)=\frac{1}{\rho(\lambda)},$
$$
\mathbb{G}_\lambda(o, \, o\, |\, z)
=
\frac{2(d-1)}{d-1+\lambda}\frac{1}{\frac{d-1-\lambda}{d-1+\lambda}+\sqrt{1-\rho(\lambda)^2z^2}}
=
\frac{a(\lambda)}{b(\lambda)+\sqrt{1-\rho(\lambda)^2z^2}}.
$$
noindent Let
$$
\Phi(t)
:=
\Phi_\lambda(t)
=
\frac{-a(\lambda)b(\lambda)+\sqrt{a(\lambda)^2+\rho(\lambda)^2(1-b(\lambda)^2)t^2}}{1-b(\lambda)^2},
\qquad
t\in\mathbb{R}.
$$
\noindent Then for any $|z|\le R_{\mathbb{G}}(\lambda)$,
$$
\mathbb{G}_\lambda(o, \, o\, |\, z)
=
\Phi\left(z \, \mathbb{G}_\lambda(o, \, o\, |\, z)\right).
$$

\noindent Define
$$
\Psi(u,\, v)
:=
\Phi(uv)-v,
\qquad
u, \, v\in\mathbb{R}\, .
$$
\noindent Then
\begin{eqnarray*}
 && \frac{\partial \Psi(u,\, v)}{\partial v}
    \Big|_{(u,\, v)=\left( \frac{1}{\rho(\lambda)}, \, \mathbb{G}_{\lambda} (o, \, o \, | \, \frac{1}{\rho(\lambda)}) \right)}
    = 0,
    \\
    \,
    \\
 &&c_1(\lambda)
    :=
    \frac{\partial^2\Psi(u, \, v)}{\partial v^2}
    \Big|_{(u,\, v)=\left( \frac{1}{\rho(\lambda)}, \, \mathbb{G}_{\lambda} (o, \, o \, | \, \frac{1}{\rho(\lambda)}) \right) }
    =
    \frac{(d-1-\lambda)^3}{2(d-1)(d-1+\lambda)^2} \not= 0,
    \\
    \,
    \\
&&c_2(\lambda)
    :=
    \frac{\partial\Psi(u, \, v)}{\partial u}
    \Big|_{(u,\, v)=\left( \frac{1}{\rho(\lambda)}, \, \mathbb{G}_{\lambda} (o, \, o \, | \, \frac{1}{\rho(\lambda)}) \right) }
    =
    \frac{2\rho(\lambda)(d-1)}{d-1-\lambda}\not= 0.
\end{eqnarray*}

\noindent Applying the method of Darboux (see \cite{BE1974} Theorem 5), we obtain that
$$
p^{(2n)}_\lambda(o,\, o)
\sim
\left(\frac{c_1(\lambda)}{2\pi\rho(\lambda)c_2(\lambda)}\right)^{1/2}\rho(\lambda)^{2n}{(2n)}^{-3/2}
=
\frac{(d-1-\lambda)^2}{16(\pi\lambda)^{1/2}(d-1)^{3/2}}\rho(\lambda)^{2n}n^{-3/2}.
$$

\noindent The idea of using the method of Darboux to establish the asymptotics for $p^{(2n)}_\lambda(o,\, o)$ is not new. For example, in Woess~\cite{WW2000} Chapter III Section 17 pp.~181--189, examples of random walk on groups are given such that $p^{(n)}(o,\, o)\sim c\rho^{n}n^{-3/2}$ for some constant $c>0$. The exact value of $c$ is not known in general.

For $z\in (-1,\, 1),$ $\mathbb{G}_{d-1}(o, \, o\, |\, z) =\frac{1}{\sqrt{1-z^2}}
=\sum_{n=0}^{\infty} \frac{(2n)!}{2^{2n}(n!)^2}z^{2n}$. Thus
$$
p_{d-1}^{(2n)}(o,\, o)
=
\frac{(2n)!}{2^{2n}(n!)^2}
\sim
\frac{1}{\sqrt{\pi n}}.
$$

\bigskip
\noindent\textbf{Acknowledgements.} The authors would like to thank an anonymous referee and Zhan Shi, Kainan Xiang, Longmin Wang for valuable comments and suggestions to improve the quality of the paper.HS's research is supported by Jiangsu University Natural Science Research Project 21KJB110003 and by Xiang Yu Ying Cai 31SH002.

\bibliographystyle{unsrt}
\bibliography{Refp5.bib}%

\begin{thebibliography}{10}

\bibitem{LR-PY2016}
R.~Lyons and Y.~Peres.
\newblock Probability on trees and networks: Plate section.
\newblock 10.1017/9781316672815, 2016.

\bibitem{BA-SA-1985}
Alberto Berretti and Alan~D. Sokal.
\newblock New monte carlo method for the self-avoiding walk.
\newblock {\em Journal of Statistical Physics}, 40(3-4):483--531, 1985.

\bibitem{LG-SA1988}
G.~F. Lawler and A.~D. Sokal.
\newblock Bounds on the l2 spectrum for markov chains and markov processes: A
  generalization of cheeger's inequality.
\newblock {\em Transactions of the American Mathematical Society},
  309(2):557--580, 1988.

\bibitem{SA-JM1989}
Sinclair~Mark Jerrum.
\newblock Approximate counting, uniform generation and rapidly mixing markov
  chains.
\newblock {\em Information and Computation}, 1989.

\bibitem{RD1994}
D.~Randall.
\newblock Counting in lattices: Combinatorial problems from statistical
  mechanics.
\newblock {\em University of California at Berkeley}, 1998.

\bibitem{BM-DD1983}
M.~Barma and D.~Dhar.
\newblock Directed diffusion in a percolation network.
\newblock {\em Journal of Physics C Solid State Physics}, 16(8):1451, 2000.

\bibitem{DD1984}
D.~Dhar.
\newblock Diffusion and drift on percolation networks in an external field.
\newblock {\em Journal of Physics A General Physics}, 17(5):L257, 1984.

\bibitem{DD-SD1998}
Stauffer~D. Dhar~D.
\newblock Drifit and trapping in biased diffusion on disordered lattices.
\newblock {\em IInternational Journal Of Modern Physics C}, 9(2):349--355,
  1998.

\bibitem{HS-BA2002}
Shlomo Havlin and Daniel Ben-Avraham.
\newblock Diffusion in disordered media.
\newblock {\em Chemometrics Intelligent Laboratory Systems},
  10(1–2):117--122, 1988.

\bibitem{SSSWX2021}
H.~Song L. Wang K. N.~Xiang Z.~Shi, V.~Sidoravicius.
\newblock Uniform spanning forests associated with biased random walks on
  euclidean lattices.
\newblock {\em ANN I H POINCARE-PR}, 57(3):1569--1582, 2021.

\bibitem{LR1990}
R.~Lyons.
\newblock Random walks and percolation on trees.
\newblock {\em The Annals of Probability}, 18(3):931--958, 1990.

\bibitem{LR1992}
Lyons and Russell.
\newblock Random walks, capacity and percolation on trees.
\newblock {\em Annals of Probability}, 20(4):2043--2088, 1992.

\bibitem{LR1995}
Lyons and Russell.
\newblock Random walks and the growth of groups.
\newblock {\em Comptes Rendus de l'Académie des Sciences - Series I -
  Mathematics}, 320:1361--1366, 1995.

\bibitem{LR-PR-PY1996a}
Russell Lyons, Robin Pemantle, and Yuval Peres.
\newblock Random walks on the lamplighter group.
\newblock {\em Annals of Probability}, 24(4):1993--2006, 1996.

\bibitem{LR-PR-PY1996b}
Russell Lyons, Robin Pemantle, and Yuval Peres.
\newblock Biased random walks on galton-watson trees.
\newblock {\em Probab Theory Rel}, 106(2):249--264, 1996.

\bibitem{BG-HY-OS2013}
G.~B. Arous, Y.~Hu, S.~Olla, and O.~Zeitouni.
\newblock Einstein relation for biased random walk on galton--watson trees.
\newblock {\em ANN I H POINCARE-PR}, 2013,49(3)(-):698--721, 2013.

\bibitem{AE2014}
Aidekon and E.
\newblock Speed of the biased random walk on a galton-watson tree.
\newblock {\em Probab Theory Rel}, 2014.

\bibitem{BG-FA-SV2014}
Gérard Ben~Arous, Alexander Fribergh, and Vladas Sidoravicius.
\newblock Lyons‐pemantle‐peres monotonicity problem for high biases.
\newblock {\em Communications on Pure Applied Mathematics}, 67(4):519--530,
  2014.

\bibitem{HY-SZ2015}
YY and Shi.
\newblock The most visited sites of biased random walks on trees.
\newblock {\em ELECTRON J PROBAB}, 2015.

\bibitem{BG-FA2014}
G.~B. Arous and A.~Fribergh.
\newblock Biased random walks on random graphs.
\newblock {\em Lecture Notes of the Institute for Computer Sciences Social
  Informatics Telecommunications Engineering}, 3(4):95--106, 2014.

\bibitem{BE1974}
E.~A. Bender.
\newblock Asymptotic methods in enumeration.
\newblock {\em Siam Review}, 16(4):485--515, 1974.

\bibitem{WW2000}
W.~Woess.
\newblock {\em Random walks on infinite graphs and groups}.
\newblock 2000.

\end{thebibliography}

\end{document}